\documentclass[11pt]{amsart}
\usepackage[normalem]{ulem}
\usepackage[usenames]{color}
\usepackage{graphicx,amscd,psfrag,amssymb,pdfpages}
\usepackage[colorlinks=true,citecolor=red,linkcolor=blue]{hyperref}

\begin{document}
\begin{center}
\textbf{\LARGE Gabriel's Paper Horn}\\
David Richeson\\   
Dickinson College\\ 
richesod@dickinson.edu
\end{center}
\vspace{.25in}

\begin{figure}[ht]
\includegraphics[width=4in]{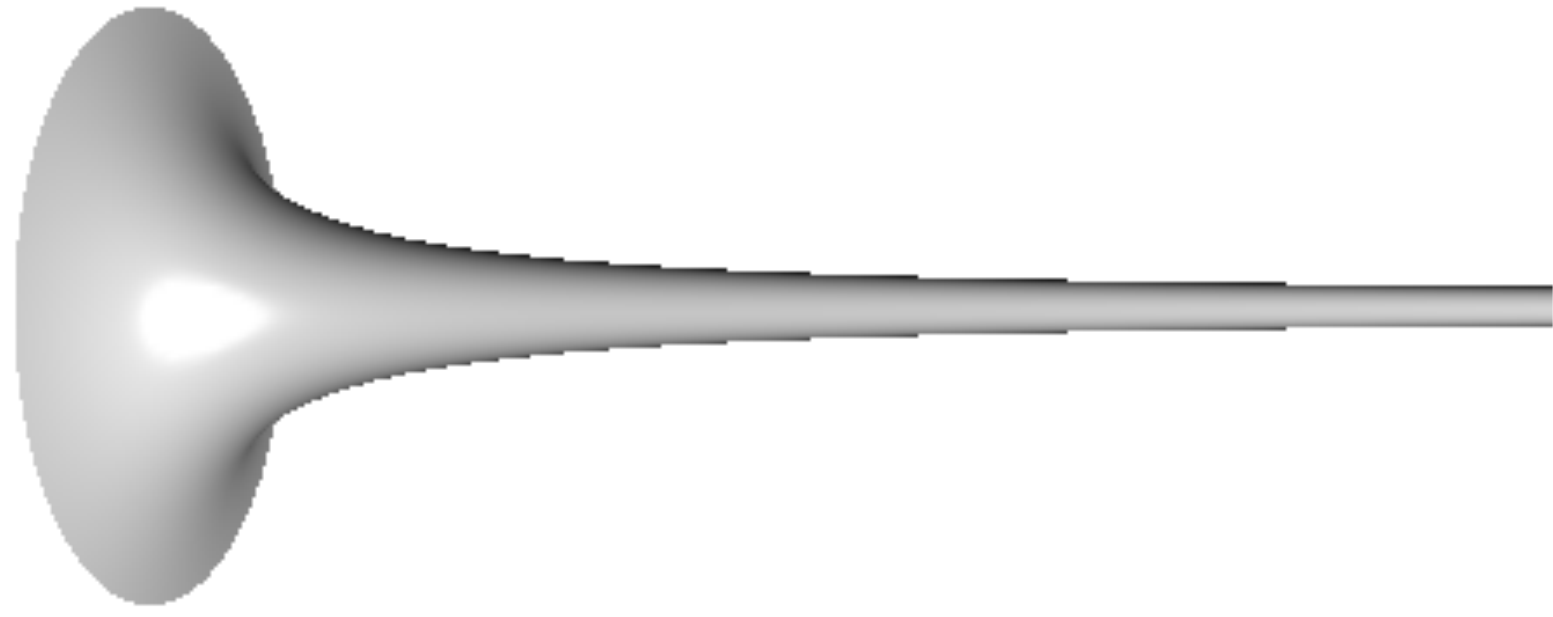}
\caption{Gabriel's horn}
\label{fig:gabrielshorn}
\end{figure}

\emph{Gabriel's horn} is the surface obtained by revolving the curve $y=1/x$ ($x\ge 1/2$) about the $x$-axis (see figure \ref{fig:gabrielshorn}). Mathematics professors wow introductory calculus students by sharing its paradoxical properties: it has finite volume, but infinite surface area. As they say, ``you can fill it with paint, but you can't paint it.''

\begin{figure}[ht]
\includegraphics[height=2.75in]{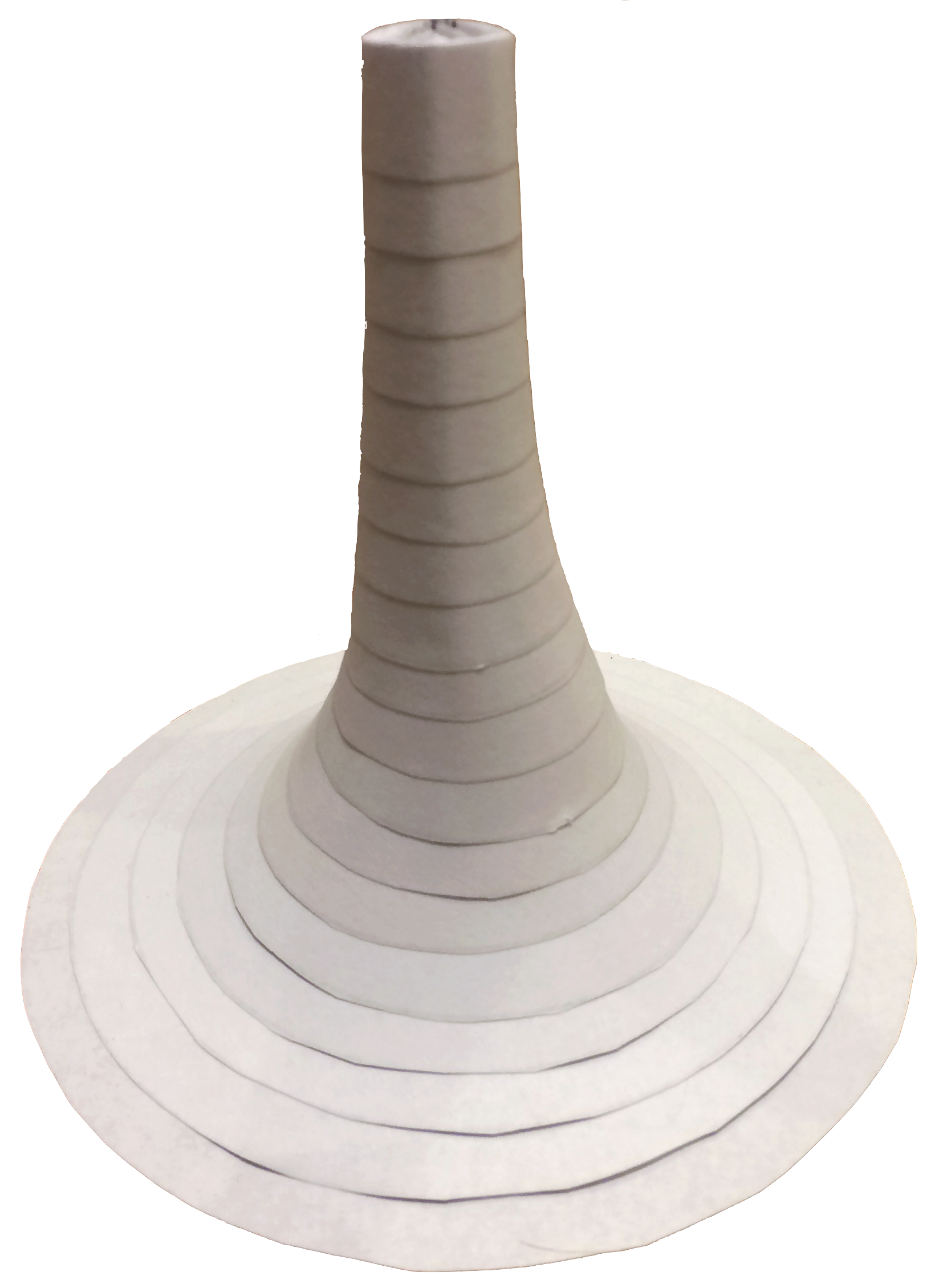}\,\,\,\,
\includegraphics[height=2.75in]{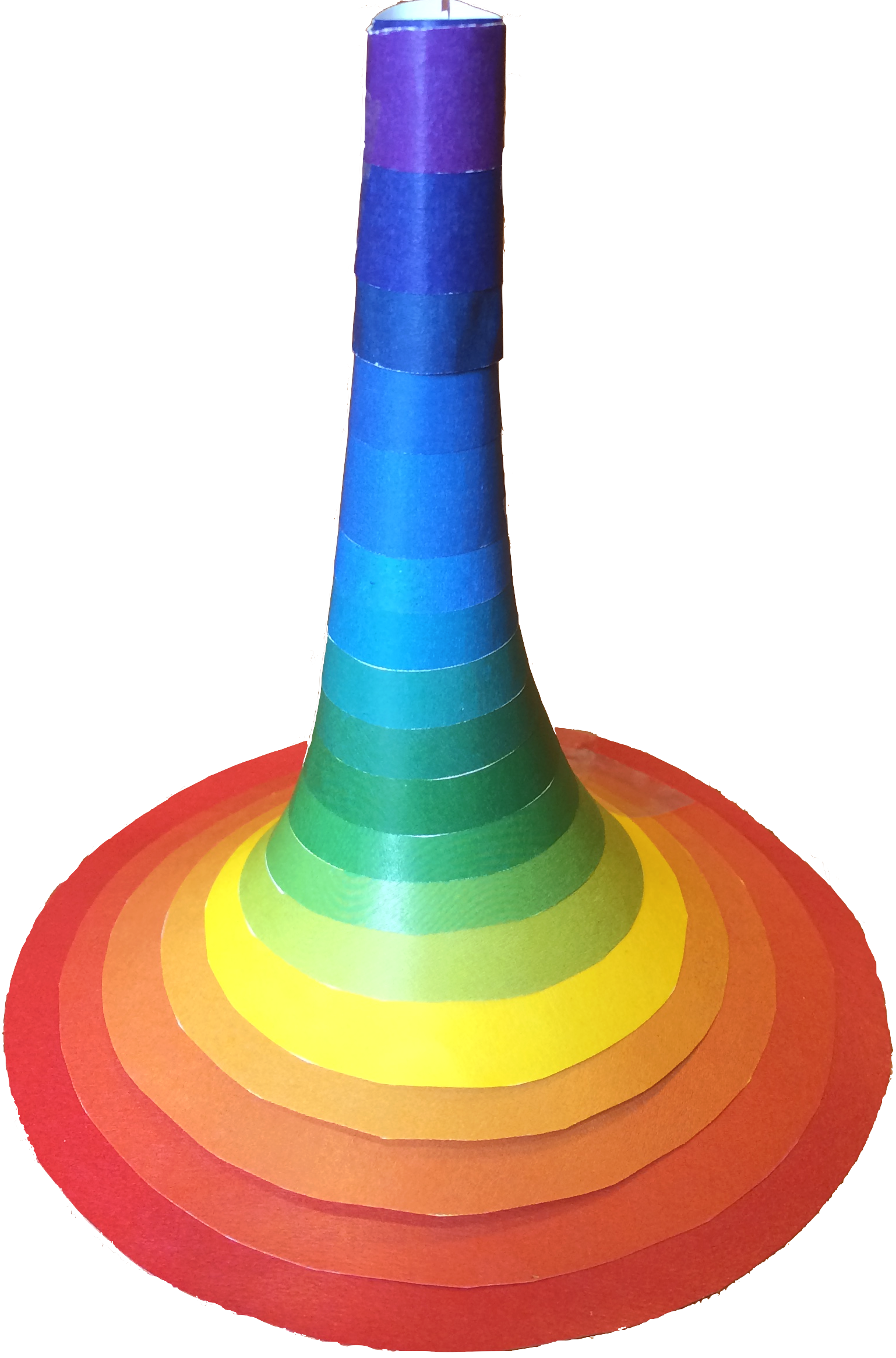}\\
\includegraphics[height=2.75in]{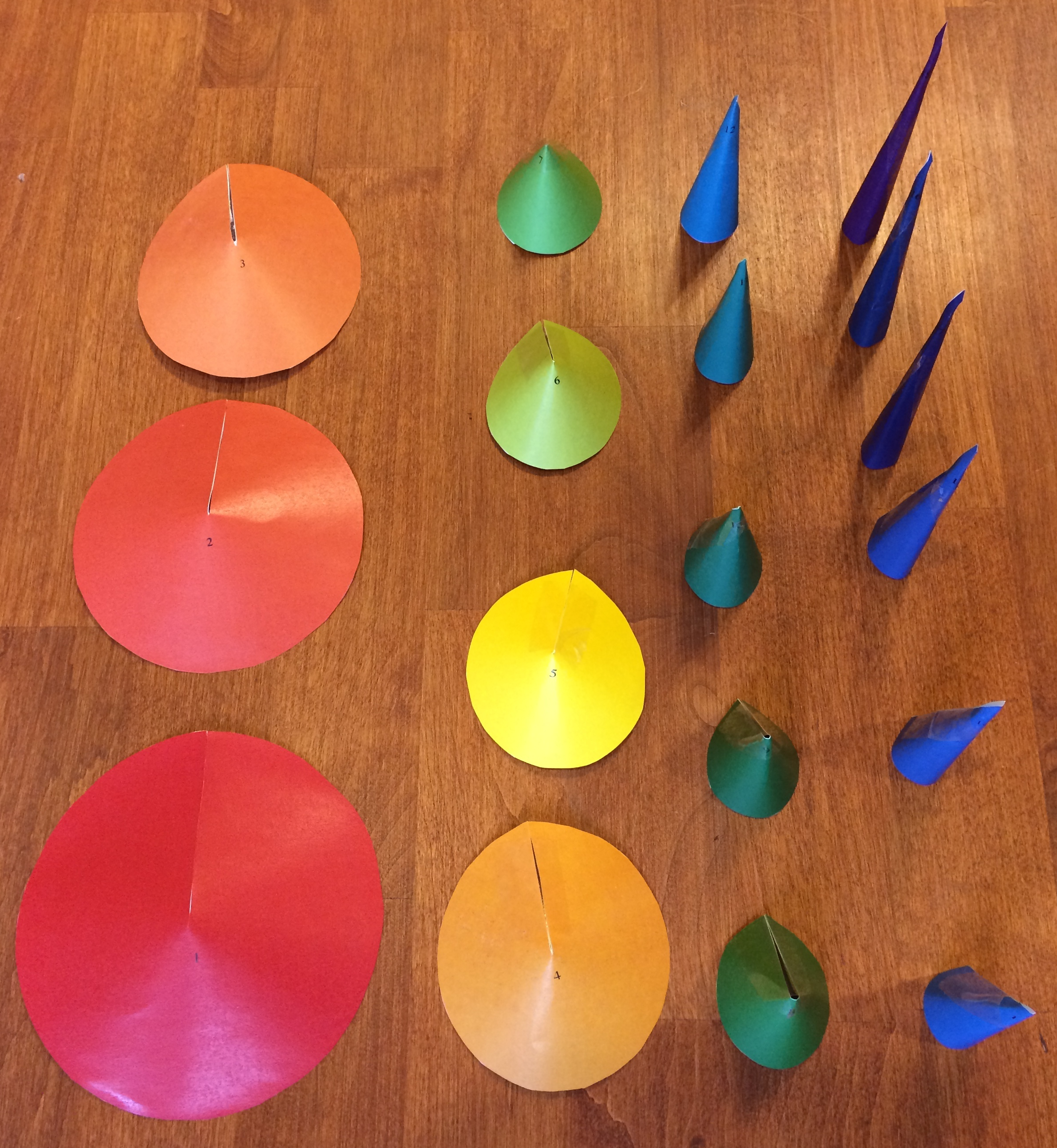}
\caption{Gabriel's horn made from paper cones}
\label{fig:hornphoto}
\end{figure}

At the end of this article we provide two templates (one color, one white) for making a model of Gabriel's horn out of paper cones, such as the one in figure \ref{fig:hornphoto}. The instructions are simple: Cut out the sectors, tape them together to form cones, and stack the cones in numerical order. The approximation is good only for the first 1/4'' of the last cone, so the stack of cones should be cut off at that height.\\

\begin{figure}[ht]
\includegraphics[width=4in]{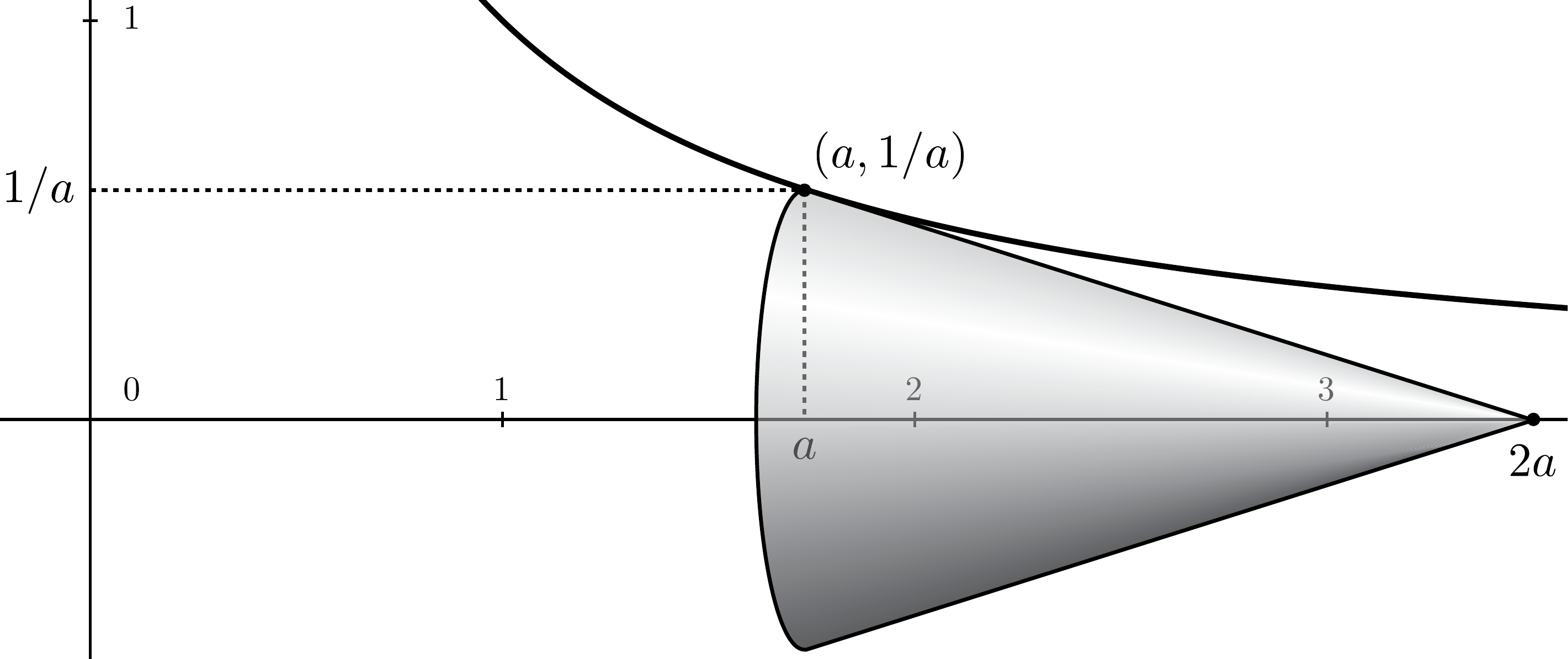}
\caption{A tangent line segment revolved to obtain a cone}
\label{fig:oneline}
\end{figure}

The rest of this article describes the mathematics used to make the sectors. In short: Take evenly spaced points along the curve, find the segments of the tangent lines between these points and the $x$-axis, and use them to generate the cones (see figure  \ref{fig:oneline}). We use $y=1/x$ as the generating curve, but this procedure works for any curve that is positive, decreasing, and concave up.

First we use calculus to find the equation of the tangent line to the graph at the point $(a,1/a)$: \[y-\frac{1}{a}=-\frac{1}{a^{2}}(x-a).\] Then we set $y=0$ to find the $x$-intercept of the tangent line: $x=2a$. From this we conclude that the radius of the base of the corresponding cone is $1/a$ and the slant height is $\sqrt{a^{2}-1/a^{2}}$. 

Now, imagine cutting the cone and unrolling it into a sector with central angle $\theta$. The arc of the sector is the circumference of the base of the cone, $2\pi/a$. The radius of the sector is the slant height of the cone, so the circumference of a full paper disk with this radius is $2\pi\sqrt{a^{2}+(1/a)^{2}}$. We use ratios to find $\theta$: 
\[\frac{\theta}{360^{\circ}}=\frac{\text{circumference(cone)}}{\text{circumference(paper)}}=\frac{2\pi/a}{2\pi\sqrt{a^{2}-1/a^{2}}}.\]
So, $\theta=\left(360/\sqrt{a^{4}+1}\right)^{\circ}.$

\begin{figure}[ht]
\includegraphics[width=5in]{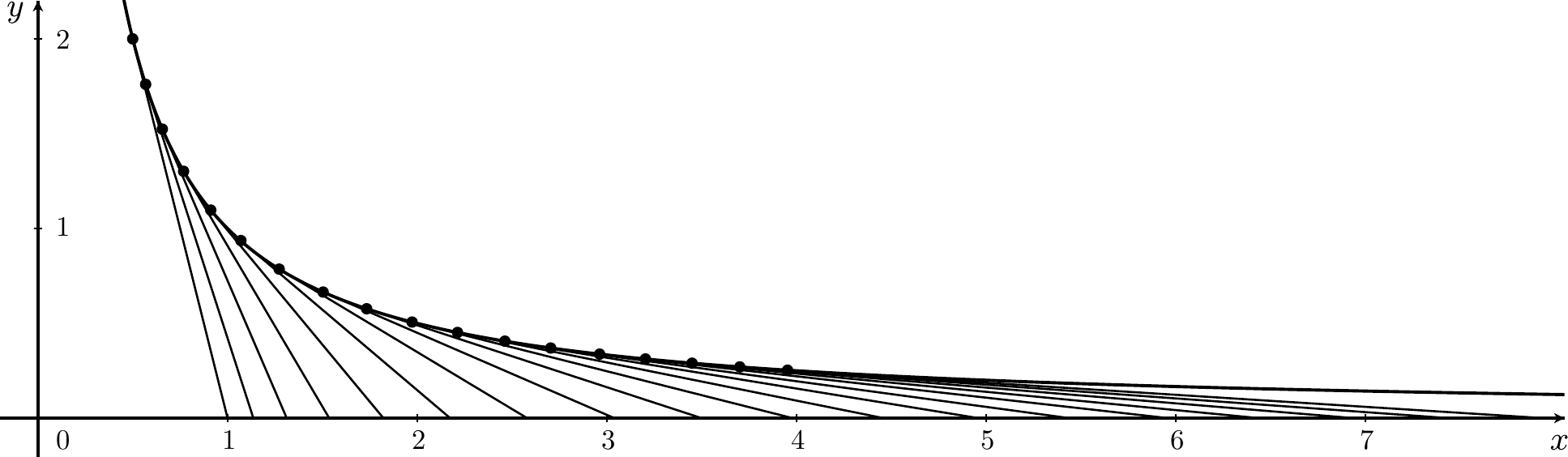}
\caption{Evenly spaced line segments tangent to $y=1/x$}
\label{fig:tanlines}
\end{figure}

We want the cones to be evenly spaced along the surface. That is, we want the visible bands to have the same widths. To accomplish this we must find points that are evenly spaced along the curve (see figure \ref{fig:tanlines}). Again we turn to calculus. The first cone corresponds to the point $(0.5,2)$. The length of the curve from $(0.5,2)$ to $(a,1/a)$ is 
\[\int_{0.5}^{a}\sqrt{1+\frac{1}{x^{4}}}\, dx.\]
We used a computer algebra system to find the seventeen $a$-values so that this integral equals $0.25, 0.5, 0.75,\ldots,4.25.$

Finally, for each of these $a$-values, we construct a paper sector with radius $\sqrt{a^{2}-1/a^{2}}$ and central angle $\left(360/\sqrt{a^{4}+1}\right)^{\circ}$.\\

This project was inspired by Daniel Walsh \cite{walsh} and Burkard Polster \cite{polster} who made pseudospheres out of paper cones. The pseudosphere is generated by revolving the tractrix about the $x$-axis. That model has the special property that the radii of the paper disks are all the same.

\includepdf{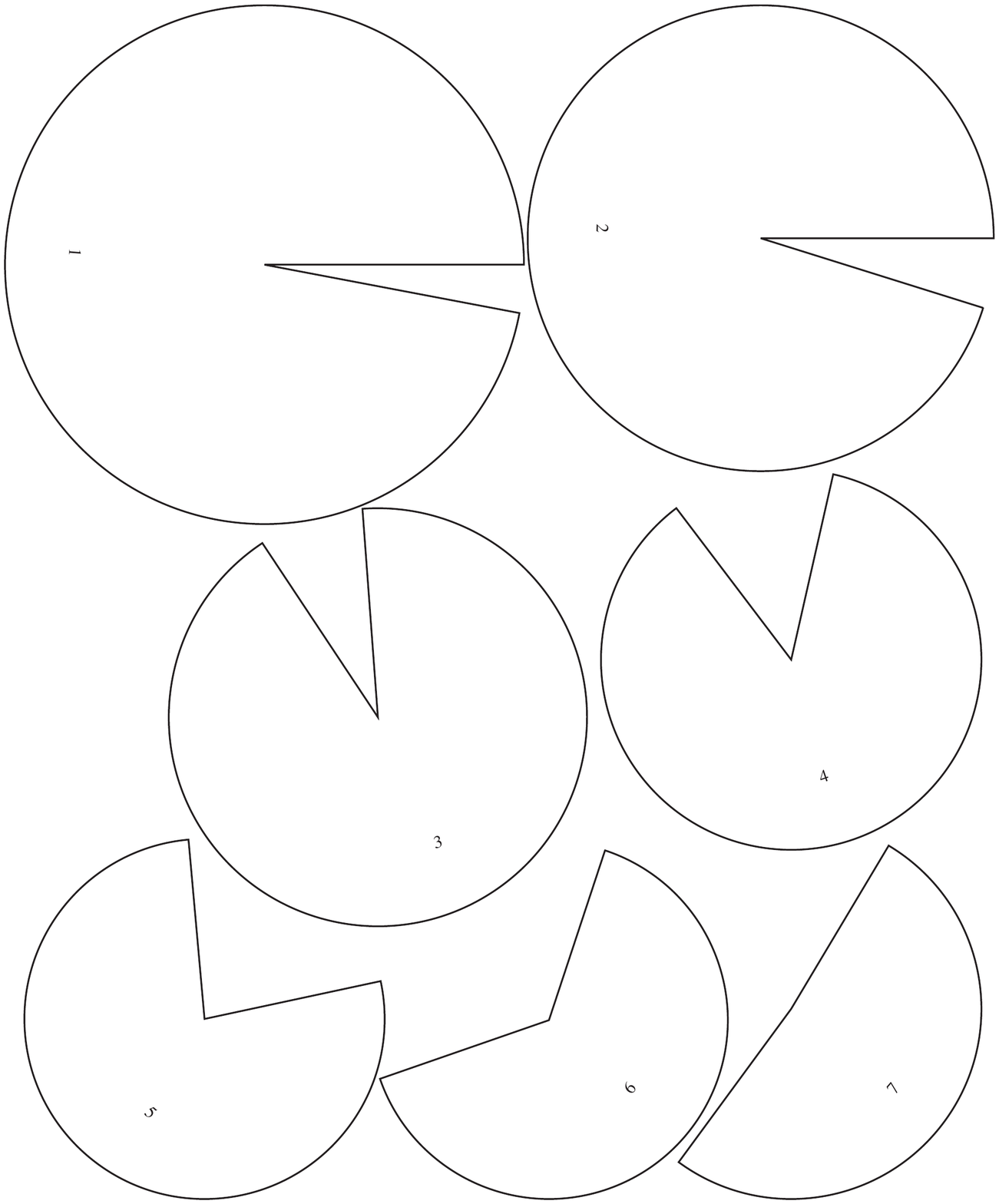}
\includepdf{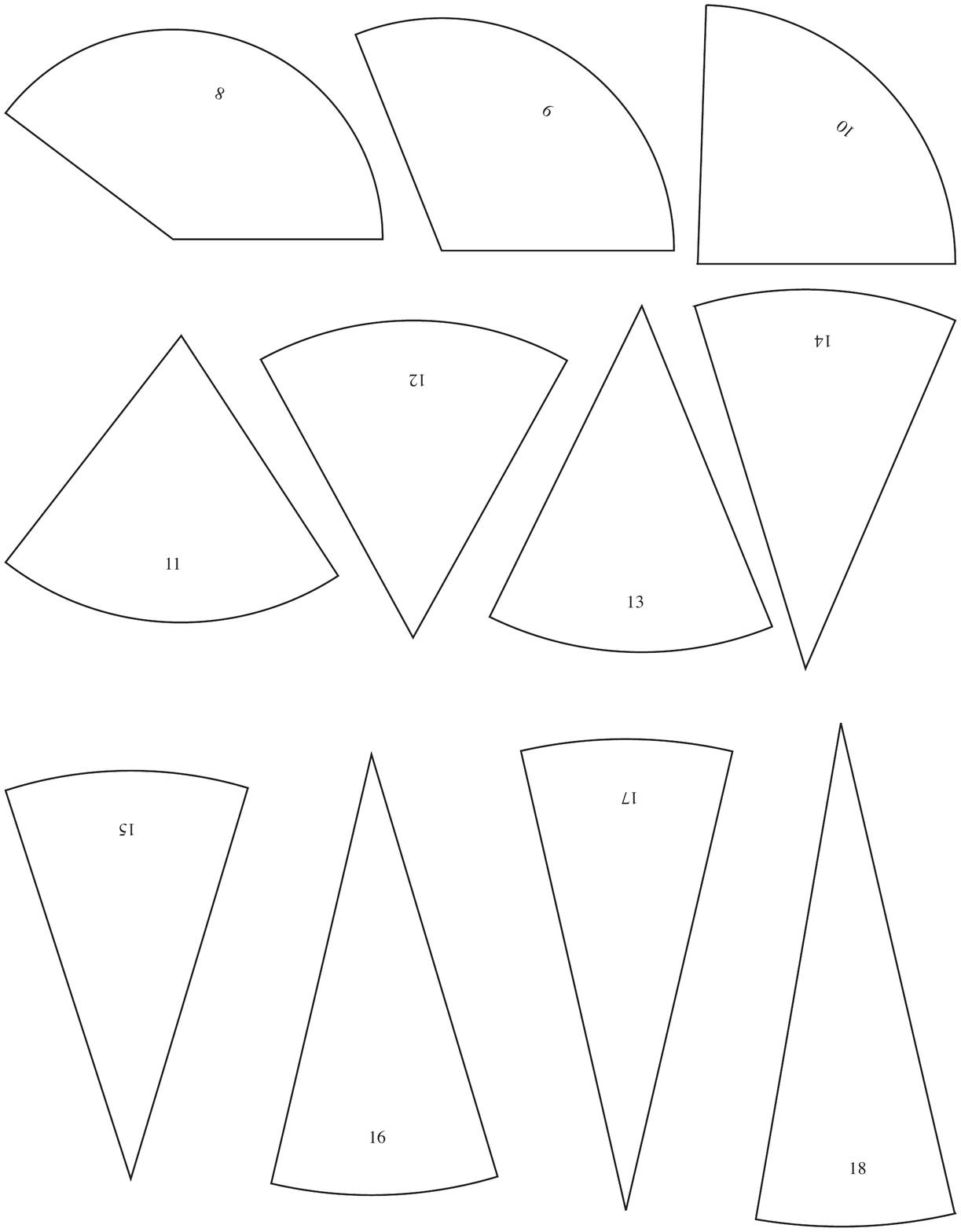}
\includepdf{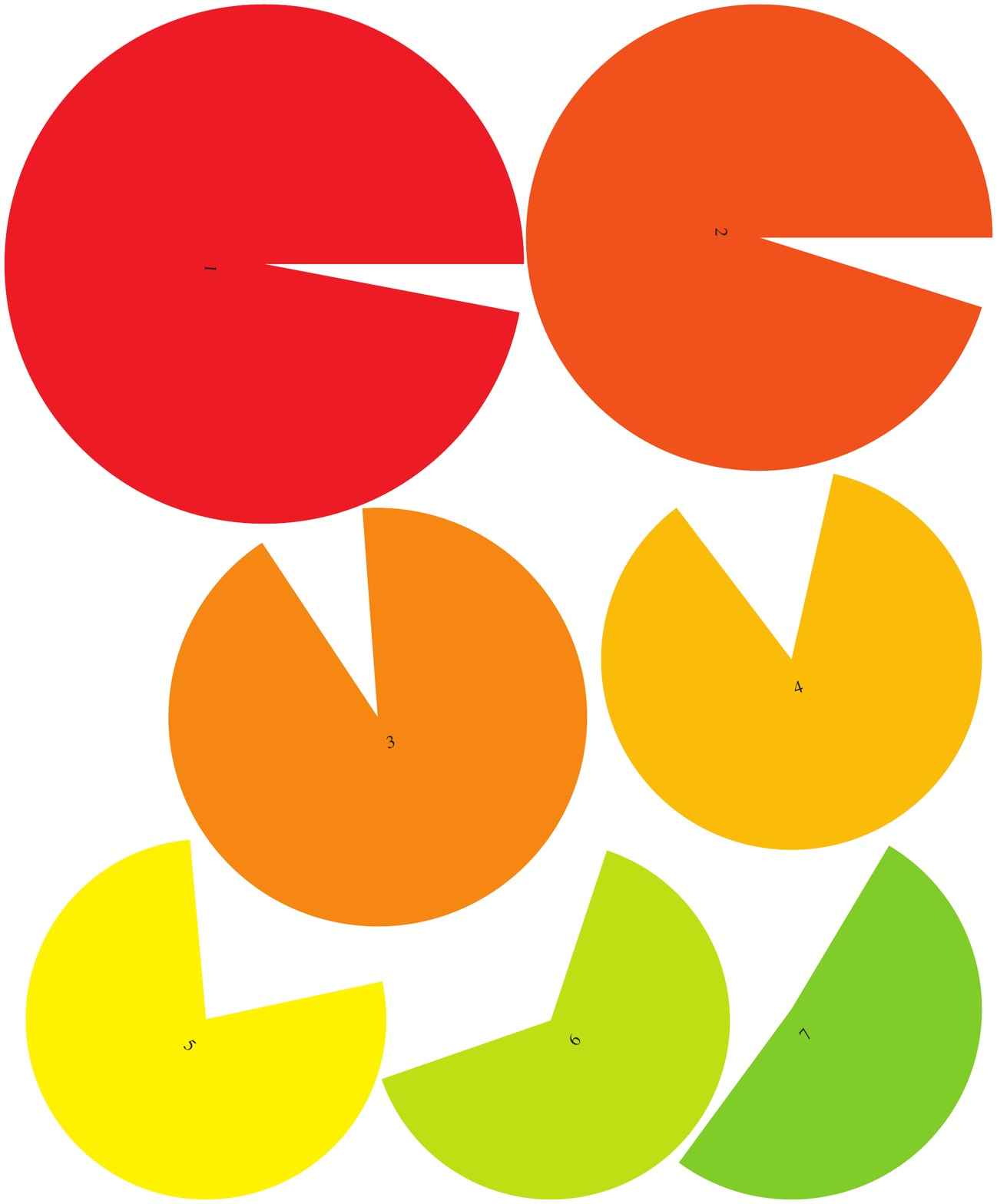}
\includepdf{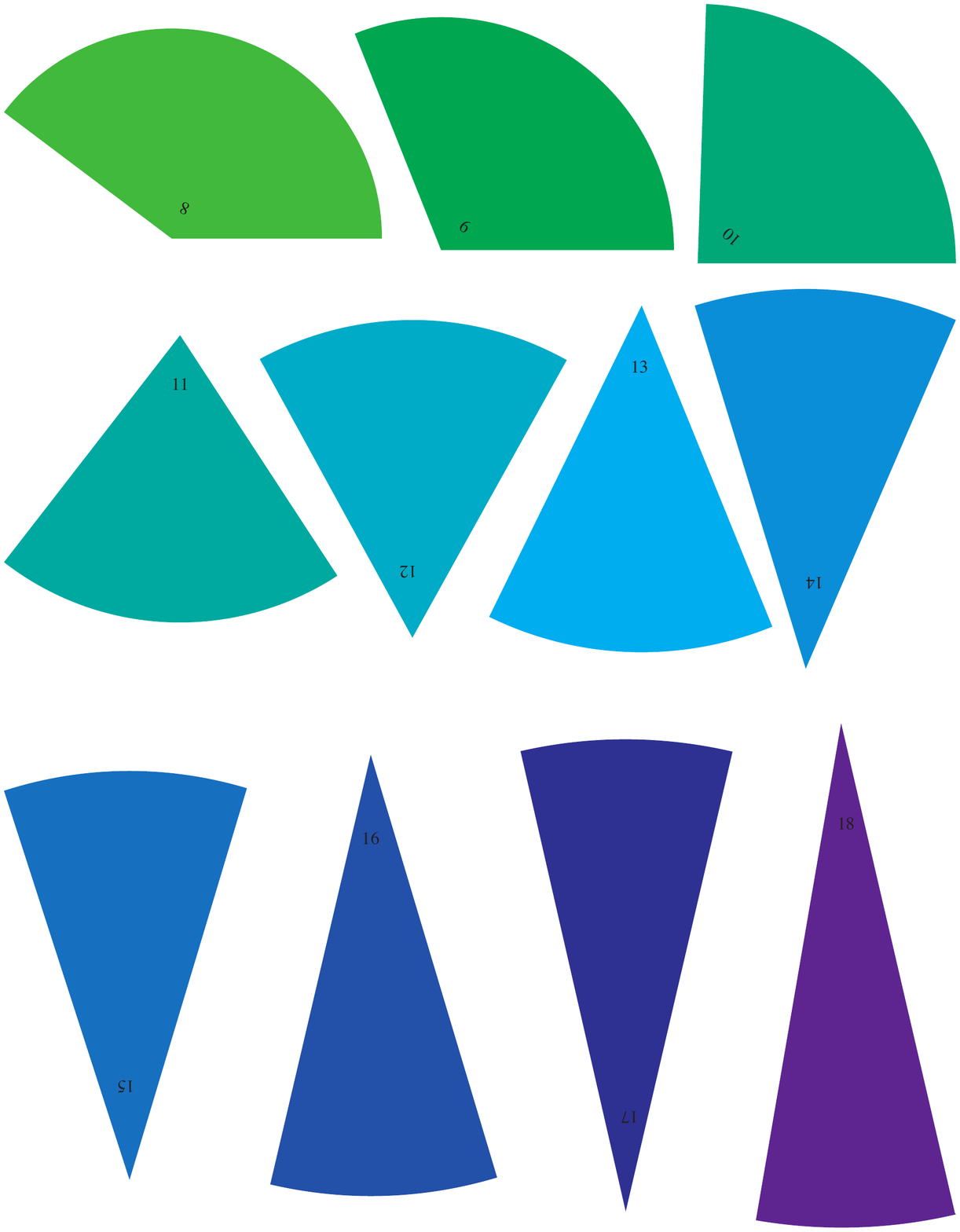}


\begin{thebibliography}{9}
\bibitem{polster} Polster, Burkard ``TracTricks,'' \emph{Math Horizons}, April 2014, pp. 18--19.
\bibitem{walsh} Walsh, Daniel, ``Sudo make me a pseudosphere,'' December 11, 2012, http://danielwalsh.tumblr.com/post/2173134224/sudo-make-me-a-pseudosphere
\end{thebibliography}
\end{document}